\documentclass[a4paper]{article}

\usepackage[english]{babel}
\usepackage[utf8x]{inputenc}
\usepackage[T1]{fontenc}

\usepackage{amssymb}
\usepackage{amsmath}
\usepackage{exscale}

\usepackage[a4paper,top=3cm,bottom=2cm,left=3cm,right=3cm,marginparwidth=1.75cm]{geometry}

\usepackage{amsmath}
\usepackage{graphicx}
\usepackage[colorinlistoftodos]{todonotes}
\usepackage[colorlinks=true, allcolors=blue]{hyperref}

\newtheorem{theorem}{Theorem}
\newtheorem{corollary}{Corollary}

\newcommand{\bA}{\mathbf{A}}

\newcommand{\bP}{\mathbf{P}}
\newcommand{\bR}{\mathbf{R}}

\newcommand{\sD}{\mathbb{D}}

\newcommand{\sR}{\mathbb{R}}

\newcommand{\ccF}{\mathcal{F}}

\newcommand{\newprocess}[1]{{#1 =(#1_t)_{0 \le t \le T}}}
\newcommand{\running}[3]{#1=#2,\ldots,#3}
\newcommand{\alli}{\running{i}{0}{n}}
\newcommand{\thetavector}{{\boldsymbol\theta}}
\newcommand{\thetaestimate}{\widehat{\thetavector}}
\newcommand{\thetaestimateml}{\widehat{\thetavector}_{{\rm ML}}}
\newcommand{\thetaestimatebayes}{\widehat{\thetavector}_{{\rm BAYES}}}

\newcommand{\taubayes}{{\tau}_{{\rm BAYES}}}
\newcommand{\deltabayes}{{\delta}_{{\rm BAYES}}}
\newcommand{\lprocess}{\Lambda^H}
\newcommand{\qfunc}{Q_H}
\newcommand{\Mprocess}{M^H}
\newcommand{\mprocess}{m^H}
\newcommand{\psiprocessn}{{\boldsymbol\psi}^H}
\newcommand{\Rfunc}{\bR_H}
\newcommand{\Zprocess}{Z^H}
\newcommand{\Cprocess}{C^H}
\newcommand{\wfunc}{w_H}
\newcommand{\wfuncdiff}{dw^H}
\newcommand{\fbm}{B^H}
\newcommand{\filtration}{\ccF^{\xi}}

\title{Optimal estimation of~a~signal
perturbed~by a~fractional Brownian noise}
\author{Artemov~A.\,V.
\and Burnaev~E.\,V.}
\date{12 February 2015}

\begin{document}
\maketitle

\begin{abstract}
We consider the problem of optimal estimation of the value of
a vector parameter~$\thetavector=(\theta_0,\ldots,\theta_n)^{\top}$
of a drift term in a fractional Brownian motion represented by a finite
sum~$\sum_{i=0}^{n}\theta_{i}\varphi_{i}(t)$ over 
known functions~$\varphi_i(t)$, $\alli$. For the value of the parameter~$\thetavector$,
we obtain a maximum likelihood estimate as well as Bayesian estimates
for normal and uniform prior distrubitions.

\end{abstract}

\section{Problem definition}\
Let~$\newprocess{\xi}$ be a stochastic process defined
on a filtered probability space 
$(\Omega,{\ccF},({\ccF}_{t})_{t\ge0},{\bf P})$
and represented by
\begin{equation}
\xi_t = a(t) + \sigma(t) \fbm_t, \label{eq:observation_process_equation}
\end{equation}
where $\newprocess{\fbm}$ is a fractional Brownian motion
with Hurst index~$H\in(0,1)$, and let drift~$a(t)$ and diffusion~$\sigma(t)$
coefficients satisfy the conditions $\int_0^T |a(t)|\,dt < \infty$ and
$\int_0^T |\sigma(t)|^{2}\,dt < \infty$, respectively. The function~$\sigma(t)$
is assumed to be known. Let the drift term~$a(t)$ be represented by a sum
\begin{equation}
\label{eq:drift_coefficient}
a(t)=\sum _{i=0}^{n}\theta_i\varphi_i(t)
\end{equation}
over the known functions~$\varphi_i(t)$, satisfying
$\int_0^T |\varphi_i(t)|\,dt < \infty$, $\alli$, with 
unknown parameters $\theta_i$, $\alli$. 
For brevity we consider vector-valued variables
$\thetavector=(\theta_0,\ldots,\theta_n)^{\top}$ and 
${\boldsymbol\varphi}(t)=(\varphi_{0}(t),\ldots,\varphi_{n}(t))^{\top}$,
such that
\begin{equation}
\label{eq:drift_coefficient_vector}
a(t)=\thetavector^{\top}{{\boldsymbol\varphi}}(t).
\end{equation}

We consider the problem of finding a sequential estimate of~$\thetavector$
given observations $\{\xi_{s},0\leqslant s\leqslant t\}$ available up to time~$t$
using the maximum likelihood and the Bayesian approaches. Within
the maximum likelihood approach, $\thetavector$ is considered as an unknown
nonrandom vector-valued parameter, and we seek to find an
estimate~$\thetaestimateml = \thetaestimateml(t)$ maximizing
the likelihood of the observed process.

In the Bayesian case, we assume~$\thetavector$ to be a random vector
taking values in~${\sR}^{n+1}$ according to some known prior distribution
$p^{\thetavector}({\bf x}), {\bf x} \in {\sR}^{n+1}$. We then
consider the problem of finding a sequential estimation rule
$\deltabayes=(\taubayes, \thetaestimatebayes)$ such that
\begin{equation}
\label{eq:optimality_criteria}
\inf_{\delta\in{\sD}}
{\bf E}\,\left[c\tau + \Vert \thetavector - \thetaestimate\Vert ^{2}\right] = {\bf E}\,\left[c\taubayes + \Vert \thetavector - \thetaestimatebayes\Vert^{2}\right],
\end{equation}
where ${\sD} = \{\delta:\,\delta =(\tau,\thetaestimate)\}$ is
a class of stopping rules with finite stopping times
$\tau\le T<\infty$ w.r.t. filtration $\filtration_t = \sigma(\{ \xi_{s},0\le s\le t\} )$.
The constant~$c > 0$ is interpreted as a cost of the observations.
The Bayesian estimation strategy consists in stopping the observations
at a time~$\taubayes$ and declaring~$\thetaestimatebayes$ to be
the optimal estimate of~$\thetavector$.

The problem of extracting a deterministic signal from observations perturbed
by a fractional Gaussian noise has attracted little attention in literature 
devoted to optimal estimation. The only work in this direction known to us
is~\cite{ccetin2013bayesian} where an optimal Bayesian estimate
for the parameter~$\mu$ of the fractional Bayesian motion with
a linear drift~$a(t) = \mu t$ is derived assuming that~$\mu$
is a normally distributed random variable with known mean and variance.

\section{Fractional Brownian motion}\
The process of fractional Brownian motion (FBM) was introduced
by Kolmogorov~\cite{kolmogorov1940spiral}
and later constructively defined by Mandelbrot~\cite{mandelbrot1968fractional}.
We use notations from~\cite{kleptsyna2000parameter}.

The standard fractional Brownian motion~$\newprocess{\fbm}$
on $[0,T]$ with Hurst index $H\in(0,1)$ is a Gaussian process
with continuous sample path such that
$$
\fbm_0 = 0,
\qquad
{\bf E}\,\fbm_t = 0,
\qquad
{\bf E}\,\fbm_s \fbm_t = \frac{1}{2}\left(s^{2H}-t^{2H}+|t-s|^{2H}\right).
$$
When $H={1}/{2}$, FBM reduces to an ordinary Brownian motion, however,
when $H\neq {1}/{2}$, FBM is not a martingale. Let us denote for $0 \le s<t\le T$
\begin{eqnarray}
\kappa_{H} & = & 2H\Gamma\bigg(\frac{3}{2}-H\bigg)\Gamma\bigg(\frac{1}{2}+H\bigg),\quad k_{H}(t,s)
 =  \kappa_{H}^{-1}s^{1/2-H}(t-s)^{1/2-H}, \nonumber\\
\lambda_{H} & =&  \frac{2H\Gamma(3-2H)\Gamma({1}/{2}+H)}{\Gamma({3}/{2}-H)},\quad\wfunc(t)
  =  \lambda_{H}^{-1}t^{2-2H},\label{eq:fbm_notation} \\
\wfuncdiff_t & =& d(\wfunc(t)) = \lambda_{H}^{-1}(2-2H)t^{1-2H}\,dt, \nonumber
\end{eqnarray}
and define the process~$\newprocess{\Mprocess}$ according to the relation
\begin{equation}
\label{eqeq23}
\Mprocess_t \equiv \int _{0}^{t}k_{H}\left(t,s\right)\,d\fbm_s.
\end{equation}
The process $\Mprocess_t$ defined in this way is a Gaussian martingale
and has the quadratic variation~$\left\langle \Mprocess_t \right\rangle$
equal to $\wfunc(t)$~(see~\cite{kleptsyna2000parameter}, \cite{norros1999elementary}).
For convenience we also define the process~$\newprocess{\mprocess}$
by the relation $\mprocess_t = \Mprocess_t / \wfunc(t)$.

\section{The Girsanov theorem for the FBM}\
In this section, we cite a result from~\cite{kleptsyna2000parameter} regarding
the likelihood process for the fractional Brownian motion with a drift.
Let~$\newprocess{Y}$ be a process defined on the filtered
probability space~$(\Omega,{\ccF},({\ccF}_{t})_{t\ge0},{\bP})$
and let its stochastic differential satisfy the relation
$$
dY_{t}=C(t)\,dt+D(t)\,d\fbm_t,
$$
where~$\newprocess{\fbm}$ is a FBM with Hurst index~$H\in(0,1)$,
and functions~$C(t)$ and~$D(t)$ are such that the function~$\qfunc(t)$
is properly defined by the relation
\begin{equation}
\label{eq:auxiliary_q}
\qfunc(t) =\frac{d}{\wfuncdiff_t}\int _{0}^{t}k_{H}(t,s)D^{-1}(s)C(s)\,ds.
\end{equation}
In this formula, differentiation w.r.t. $\wfuncdiff_t$
is understood in the following way:
$$
\frac{df(t)}{\wfuncdiff_t} \equiv  
\frac{\lambda_H}{2-2H}\,t^{2H - 1}\,\frac{df(t)}{dt}.
$$
Defining the function~$\qfunc(t)$ allows one to formulate 
an analogue of the Girsanov theorem for the process~$Y$.

\begin{theorem}[{see\ \cite{kleptsyna2000parameter}}]
Let~$\qfunc(t)$ belong to the space $L^{2}([0,T],\wfuncdiff_t)$,
where the quantity $\wfuncdiff_t$ is defined by~{\rm\eqref{eq:fbm_notation}}.
Let us define a random process~$\newprocess{\lprocess}$ by the relation
\begin{equation}
\label{eq:fractional_likelihood}
\lprocess_t = \exp\bigg\{
    \int_{0}^{t}\qfunc(s)\, d\Mprocess_s
    -\frac{1}{2}\int _{0}^{t}(\qfunc(s))^{2}\wfuncdiff_s\bigg\}.
\end{equation}
The ${\bf E} \,\lprocess_t = 1$ and the distribution of $Y$ w.r.t.
the measure ${\bf P}^Y=\lprocess_t {\bf P}$ coincides with
the distribution of $\int _{0}^{t}D(s)\,d\fbm_s$ w.r.t. ${\bf P}$.
\end{theorem}

The random process $\lprocess$ is called the likelihood process
or the Radon-Nikodym derivative ${d{\bf P}^Y}/{d{\bf P}}$
of the measure ${\bf P}^Y$ w.r.t. the measure ${\bf P}$.

\section{The maximum likelihood estimate of the drift parameter}\
Let us consider the problem of finding the maximum likelihood
estimate for the drift parameter~$\thetavector$ defined
in~\eqref{eq:observation_process_equation}.
According to~\eqref{eq:observation_process_equation}--\eqref{eq:drift_coefficient},
the process~$\xi_t$ satisfies the equation
\begin{equation}
\label{eq:observation_explicit_form}
\xi_{t}=\sum _{i=0}^{n}\theta_{i}\varphi_{i}(t)+\sigma(t)\fbm_t,
\end{equation}
while its stochastic differential satisfies the relation
$$
d\xi_{t}=\sum _{i=0}^{n}\theta_{i}\varphi_{i}'(t)\,dt+\sigma(t)\,d\fbm_t.
$$
The structure of the likelihood process and the corresponding estimate is described
by the following theorem.

\begin{theorem}
Let the drift coefficient $a(t)$ of the fractional Brownian motion
have the form~{\rm\eqref{eq:drift_coefficient}}--{\rm\eqref{eq:drift_coefficient_vector}}.
Then the maximum likelihood estimate~$\thetaestimateml$
for the drift parameter~$\thetavector$ is defined by
\begin{equation}
\label{eq:ml_estimate_generic}
\thetaestimateml = \Rfunc^{-1}(t) \psiprocessn_t,
\end{equation}
where~$\Rfunc(t)$ is a nonrandom matrix with elements defined by
\begin{equation}
\label{eq:psi_H_formula_R}
(\Rfunc(t))_{ij} =  \int _0^t\psi_i(s)\psi_j(s)\,\wfuncdiff_s,
\qquad \running{i,j}{0}{n},
\end{equation}
and $\newprocess{\psiprocessn}$ is a stochastic process taking values
in ${\sR}^{n+1}$ with coordinates defined by
\begin{equation}
\label{eq:psi_H_formula_psi}
(\psiprocessn_t)_i = \int _0^t\psi_i(s)\,d\Mprocess_s,
\qquad \running{i}{0}{n},
\end{equation}
where the functions $\psi_i(t),$ $\alli,$ are given by
\begin{equation}
\label{eq:ml_psi_i}
\psi_{i}(t)=\frac{d}{\wfuncdiff_t} \int _{0}^{t}k_{H}(t,s)\sigma^{-1}(s)\varphi_{i}'(s)\,ds,
\qquad \alli,
\end{equation}
and $\Mprocess_t$ is defined  by \eqref{eqeq23} with $\newprocess{\xi}$ instead of $\newprocess{\fbm}$.
\end{theorem}

\textbf{Proof.} The general form of~$\qfunc^{\thetavector}(t)$ function is defined
by~\eqref{eq:auxiliary_q}. Using the notation from~\eqref{eq:ml_psi_i}
for the functions~$\psi_{i}(t), \alli$, we obtain
$$
\qfunc^{\thetavector}(t) = \sum _{i=0}^{n}\theta_{i}\,
\frac{d}{\wfuncdiff_t} \int _{0}^{t}k_{H}(t,s)\sigma^{-1}(s)\varphi_{i}'(s)\,ds =\sum _{i=0}^{n}\theta_{i} \psi_{i}(t).
$$
The likelihood process~$\lprocess$ is then defined as
(see~\eqref{eq:fractional_likelihood}):
\begin{equation}
\label{eq:fractional_likelihood_series}
\lprocess_t(\thetavector)=\exp\bigg\{
    \sum _{i=0}^{n}\theta_{i}\int _{0}^{t}\psi_{i}(s)\,d\Mprocess_s
    - \frac{1}{2}\int _{0}^{t}\bigg(\sum _{i=0}^{n}\theta_{i}\psi_{i}(s)\bigg)^{2}\wfuncdiff_s \bigg\}.
\end{equation}
The process~$\lprocess$ defines the Radon-Nikodym derivative of the measure
generated by the observations~$\xi$ from~\eqref{eq:observation_explicit_form}
w.r.t. the measure of the process~$\xi_s = \fbm_s$, $s \leqslant t$.
Using the vector notation from~\eqref{eq:drift_coefficient_vector},
one can write the formula in~\eqref{eq:fractional_likelihood_series}
more compactly:
\begin{equation}
\label{eq:fractional_likelihood_series_vector}
\lprocess_t(\thetavector) =
\exp\bigg\{ \thetavector^{\top}\psiprocessn_t
            -\frac{1}{2}\,\thetavector^{\top}\Rfunc(t)\thetavector
    \bigg\},
\end{equation}
where the elements of the $n \times n$ matrix~$\Rfunc(t)$
and the components of the~$(n+1)$-dimensional process~$\psiprocessn$
are defined by~\eqref{eq:psi_H_formula_R} and~\eqref{eq:psi_H_formula_psi}, respectively.
The maximum likelihood estimate~$\thetaestimateml=\arg\max_{\thetavector}\lprocess_t(\thetavector)$
is obtained as a solution of the system of linear equations
$$
\int_{0}^{t}\psi_{i}(s)\,d\Mprocess_s-
\sum_{j=0}^{n}\theta_{j}\int_{0}^{t}\psi_{i}(s)\psi_{j}(s)\,\wfuncdiff_s=0,
\qquad \alli,
$$
which could be written in a vector form
$$
\psiprocessn_t - \Rfunc(t)\thetavector = 0.
$$
If the matrix~$\Rfunc(t)$ is invertible for every~$t \geqslant 0$, then
the solution of the system is~$\thetaestimateml$
from~\eqref{eq:ml_estimate_generic}.

\begin{corollary}[the case of a polynomial drift]
Let $\varphi_{i}(t)=t^{i},$ $\alli,$ and assume the diffusion coefficient
to be constant $\sigma(t)=\sigma$. Then the observable process has
the structure $\xi_{t}=\sum_{i=0}^{n}\theta_{i}t^{i}+\sigma \fbm_t,$
functions $\psi_i(t) = \beta_{H}(i) / \sigma t^{i-1},$ $\alli,$
whereas the components of the vector-valued stochastic
process~$\psiprocessn$ from~{\rm\eqref{eq:psi_H_formula_psi}}
and the elements of the matrix~$\Rfunc(t)$ from~{\rm\eqref{eq:psi_H_formula_R}} are defined by
$$
(\psiprocessn_t)_i = \frac{\beta_{H}(i)}{\sigma}\int _{0}^{t}s^{i-1}\,d\Mprocess_s,
\quad (\Rfunc(t))_{ij} = \frac{\alpha_H(i,j)}{\sigma^2}\,t^{i+j-2H}
$$
respectively, where
\begin{eqnarray*}
\alpha_H(i,j) &=& \lambda^{-1}_{H}\beta_{H}(i)\beta_{H}(j)\,\frac{2-2H}{i+j-2H},\\
\beta_{H}(i) &=& i \,\frac{2-2H+i-1}{2-2H}\,
\frac{\Gamma(3-2H)}{\Gamma(3-2H+i-1)}\,\frac{\Gamma(3/2-H+i-1)}{\Gamma(3/2-H)},
\qquad\running{i,j}{0}{n}.
\end{eqnarray*}
The maximum likelihood estimate~$\thetaestimateml$ is obtained as a solution
to the equation~$\psiprocessn_t - \Rfunc(t)\thetavector = 0$.
\end{corollary}

We note that for $n=1$ the observable process satisfies the stochastic
differential equation~$d\xi_{t} = \theta_{1}dt+\sigma\,dB_{t}^{H}$
and the likelihood process has the form
$\lprocess_t(\thetavector) = \exp\{\theta_{1}\sigma^{-1}\Mprocess_t - \theta_{1}^{2}\sigma^{-2}\lambda_{H}^{-1}t^{2-2H} / 2\}$.
Therefore the maximum likelihood estimate~$(\widehat{\theta}_{1})_{{\rm ML}}$
of~$\theta_{1}$ has the form
\begin{equation}
(\widehat{\theta}_{1})_{{\rm ML}}=\frac{\sigma \Mprocess_t}{\wfunc(t)}.
\end{equation}
This particular result (for $\sigma = 1$) has been obtained
in~\cite{norros1999elementary}.

\section{The Bayesian estimate of the drift parameter}\
Consider the problem of finding the Bayesian estimate
of the parameter~$\thetavector\in{\sR}^{n + 1}$ assuming that~$\thetavector$
has a prior distribution~${\bP}^{\thetavector}$ with density~$p^{\thetavector}({\bf x})$,
${\bf x} = (x_0,\ldots,x_n) \in {\sR}^{n + 1}$.

According to the generalized Bayes rule (see~\cite{ccetin2013bayesian},~\cite{liptser1974statistics}),
the conditional distribution density of~$\thetavector$ given
observations~${\ccF}_{t}^{\xi} = \sigma(\{\xi_{s},\ 0\le s\le t\})$
is represented by
\begin{eqnarray}
p^{\thetavector}({\bf x}\,|\,{\ccF}_{t}^{\xi}) &=&
\frac{d{\bf P}(\theta_{0}\le x_{0},\ldots,\theta_{n}\le x_{n}\,|\,\filtration_t)}{dx_{0}\cdots dx_{n}}\nonumber\\
&=&\frac{p^{\thetavector}({\bf x})\lprocess_t({\bf x})}{\int_{{\sR}^{n + 1}}p^{\thetavector}({\bf z})\lprocess_t({\bf z})\, d{\bf z}},
 \qquad {\bf x} \in \sR^{n + 1},
\label{eq:generalized_bayes_theorem}
\end{eqnarray}
where $\lprocess_t({\bf x})$ is the likelihood process previously described in section~4.
We further consider two special cases where the prior distribution of~$\thetavector$
is either normal or uniform.

\subsection{The case of a normal prior distribution}\
The main result of this section is presented in the following theorem.

\begin{theorem}
\label{thm:bayesian_normal}
Let $\thetavector$ be a multivariate normal random variable
with mean~${\bf m}$ and covariance matrix~${\bf\Sigma}$. Then the optimal Bayesian
estimate~$\thetaestimatebayes$ for the value of~$\thetavector$ is the posterior mean
\begin{equation}
\label{eq:bayes_normal_posterior_mean}
\thetaestimatebayes = {\bf E}\,[\thetavector\,|\,\filtration_t] =
(\Rfunc(t) + {\bf\Sigma}^{-1})^{-1}(\psiprocessn_t + {\bf\Sigma}^{-1}{\bf m}).
\end{equation}
The estimation error~ ${\bf E}\,(\Vert \thetavector - \thetaestimatebayes\Vert^2\,|\,\filtration_t)$
is defined by the trace of the posterior covariance matrix
\begin{equation}
\label{eq:bayes_normal_posterior_cov}
{\rm cov}\,[\thetavector\,|\,\filtration_t] =
(\Rfunc(t) + {\bf\Sigma}^{-1})^{-1}.
\end{equation}
\end{theorem}

\textbf{Proof.} It is well known that the optimal least squares estimate for the value
of the vector~$\thetavector$ conditioned upon the observations
history~$\left\{\xi_{s},\ 0\le s\le t\right\}$ up to the moment~$t$
is defined by the conditional expectation~${\bf E}\,[\thetavector\,|\,\filtration_t]$.
The estimation error is defined by the trace of the conditional
covariance matrix~${\rm cov}\,[\thetavector\,|\,\filtration_t]$.
In what follows we show that these quantities are easy to compute
in the normal case.

Using the formulae~\eqref{eq:fractional_likelihood_series_vector}
and~\eqref{eq:generalized_bayes_theorem} and writing the 
multivariate normal density as 
$p^{\thetavector}({\bf x}) =
    (2\pi)^{-(n+1)/2}(\det{\bf\Sigma})^{-1/2}
    \exp\{ -({\bf x}-{\bf m})^{\top}{\bf\Sigma}^{-1}({\bf x}-{\bf m})/2\}$,
we obtain the following formula for the conditional distribution
of~$\thetavector$ given $\filtration_t=\sigma(\{\xi_{s},\ 0\le s\le t\})$:
$$
\begin{array}{c}
 p^{\thetavector}({\bf x}\,|\,\filtration_t) =
g({\bf x}) / \int_{{\sR}^{n + 1}}g({\bf z})\,d^{n}{\bf z},\\[10pt]
 g({\bf x}) = \exp\Big\{{\bf x}^{\top}(\psiprocessn_t + {\bf\Sigma}^{-1}{\bf m}) -{\bf x}^{\top}(\Rfunc(t) + {\bf\Sigma}^{-1}){\bf x}/2\Big\}.
\end{array}
$$
Using a well-known formula
$$
\int_{{\sR}^{n}}\exp\bigg\{
-\frac{1}{2}\,{\bf x}^{\top}\bA{\bf x}+{\bf x}^{\top}{\bf b}\bigg\}\,d^{n}{\bf x}= \sqrt{\frac{(2\pi)^{n}}{\det\bA}}\,\exp\bigg\{\frac{1}{2}\,{\bf b}^{\top}\bA^{-1}{\bf b}\bigg\}, 
$$
we rewrite the latter expression as follows:
$$
p^{\thetavector}({\bf x}\,|\,\filtration_t) =
    \sqrt{\frac{\det\bA}{(2\pi)^{n+1}}}\,
    \exp\bigg\{-\frac{1}{2}\,({\bf x}^{\top}\bA{\bf x}-2{\bf x}^{\top}{\bf b}+{\bf b}^{\top}\bA^{-1}{\bf b})\bigg\},
$$
where the nonrandom matrix~$\bA = \bA_H(t)$ and the $(n+1)$-dimensional stochastic
process $\newprocess{{\bf b}^H}$ are given by $\bA_H(t) = \Rfunc(t) + {\bf\Sigma}^{-1}$
and~${\bf b} = {\bf b}^H_t = \psiprocessn_t + {\bf\Sigma}^{-1}{\bf m}$, respectively.
It can be readily seen that the conditional density
$$
p^{\thetavector}({\bf x}\,|\,\filtration_t) =
    \sqrt{\frac{\det\bA}{(2\pi)^{n+1}}}\,
    \exp\bigg\{-\frac{1}{2}\,({\bf x}-\bA^{-1}{\bf b})^{\top}\bA({\bf x}-\bA^{-1}{\bf b})\bigg\}
$$
is multivariate normal with mean and covariance
\begin{eqnarray*}
{\bf E}\,[\thetavector\,|\,\filtration_t] & = &
    \bA^{-1}{\bf b}=
        (\Rfunc(t) + {\bf\Sigma}^{-1})^{-1}
        (\psiprocessn_t + {\bf\Sigma}^{-1}{\bf m}),\\
{\rm cov}\,[\thetavector\,|\,\filtration_t] & = &
    \bA^{-1} = (\Rfunc(t) + {\bf\Sigma}^{-1})^{-1},
\end{eqnarray*}
respectively. The quantity~${\bf E}\,(\Vert \thetavector - \thetaestimatebayes \Vert^2 \,|\,\filtration_t)$
representing the conditional mean squared estimation error in the normal case
has the form
\begin{eqnarray*}
{\bf E}\,(\Vert \thetavector - \thetaestimatebayes \Vert^2 \,|\,\filtration_t)&=&
{\bf E}\Big[{\rm tr}\,((\Rfunc(t) + {\bf\Sigma}^{-1})^{-1})\Big]\\
&=& {\rm tr}\,\left((\Rfunc(t) + {\bf\Sigma}^{-1})^{-1}\right).
\end{eqnarray*}

\begin{corollary}
\label{thm:bayesian_normal_stopping}
Let the conditions of theorem~{\rm\ref{thm:bayesian_normal}} be satisfied.
Then the optimal stopping time for~{\rm\eqref{eq:optimality_criteria}} is nonrandom.
\end{corollary}

\textbf{Proof.} In order to determine the optimal stopping time~$\taubayes$ in~{\rm\eqref{eq:optimality_criteria}}
we have to solve the following optimal stopping problem:
$$
\taubayes = \arg \inf_{\tau\in{\sD}} 
{\bf E}\,
        \left[
                c\tau + 
                {\bf E}\,\left(\Vert \thetavector - \thetaestimatebayes \Vert^{2} \,|\,{\ccF}_{\tau}^{\xi} \right)
        \right] = 
        \arg \inf_{t \in [0, T]} F_H(t),
$$
where the function
\begin{equation}
F_H(t) =
ct + {\bf E}\,(\Vert \thetavector - \thetaestimatebayes \Vert^{2} \,|\,\filtration_t ) =
ct + {\rm tr}\,\left((\Rfunc(t) + {\bf\Sigma}^{-1})^{-1}\right),\qquad t \in [0, T],
\label{eq:f_h_bayes_normal_common}
\end{equation}
is nonrandom.

\begin{corollary}[the case of a polynomial drift] 
Let $\varphi_{i}(t)=t^{i},$ $\alli,$ assume the diffusion coefficient
to be constant, $\sigma(t)=\sigma$, and let the covariance
matrix~${\bf\Sigma}$ be diagonal meaning that all $\theta_0,\ldots,\theta_n$
are independent of each other. Then the function~$F_H(t)$
from~{\rm\eqref{eq:f_h_bayes_normal_common}} has a single minimum for $t \in [0, T]$.
\end{corollary}

\textbf{Proof.} Let~${\bf\Sigma} = {\rm diag}\,(\gamma_0^2,\ldots,\gamma_n^2)$,
where~$\gamma_i^2 = {\rm var}(\theta_i), \alli$. Then the trace of the conditional
covariance matrix has the form
$$
{\rm tr}\, \left((\Rfunc(t) + {\bf\Sigma}^{-1})^{-1}\right) =
\prod_{i=0}^{n}\bigg(\frac{\alpha_H(i,i) t^{2i-2H}}{\sigma^4} + \gamma^{-2}\bigg)^{-1}
$$
and is a strictly increasing function for~$t > 0$. Thus~$F_H(t)$
in~\eqref{eq:f_h_bayes_normal_common} is a sum of a strictly
increasing and a strictly decreasing function and has a unique
minimum at some~$t \in [0, T]$.

\begin{figure}[ht]
 \centering
 \includegraphics[width=0.7\textwidth]{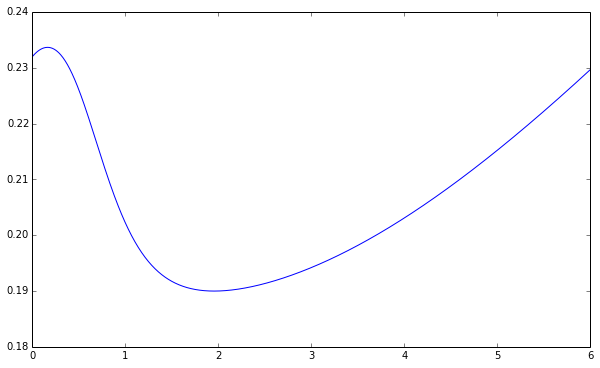}
 \caption{Graph of the cost function $F_{H}(t)$ for $n = 2,$ $H = 0.2,$ $c = 0.02$.}
 \label{fig:F_H_plot_for_n_equal_2}
\end{figure}

A similar result can be obtained for the case 
of a linear trend~$a(t) = \mu t$~\cite{ccetin2013bayesian}.
Fig. \ref{fig:F_H_plot_for_n_equal_2} presents the graph of~$F_H(t)$ for the case of a quadratic trend
and values~$H=0.2, c=0.02$.

Note that if the observable process~$\xi$ satisfies a linear
stochastic differential equation~$d\xi_{t} = \theta_{1}\,dt+\sigma \,dB_{t}^{H}$,
where~$\theta_1$ is a normally distributed random variable
with expectation~$m$ and variance~$\gamma^2$, then its conditional
distribution is normal with density
$$
p^{\theta_1}(x\,|\,\filtration_t) = 
\sqrt{\frac{\wfunc(t)/\sigma^2 + 1/\gamma^2}{2\pi}}\,
\exp\bigg\{ -\bigg(x-\frac{\Mprocess_t / \sigma + m / \gamma^2}{\wfunc(t) / \sigma^2 + 1/\gamma^2}\bigg)^2\,
\frac{\wfunc(t) / \sigma^2 + 1/\gamma^2}{2}\bigg\}.
$$
Therefore the Bayesian estimate~$(\widehat{\theta}_{1})_{{\rm BAYES}} = {\bf E}\,[\theta_{1}\,|\,\filtration_t]$
of~$\theta_{1}$ and its corresponding estimation error
${\bf E}\,[(\theta_1 - (\widehat{\theta}_{1})_{{\rm BAYES}})^2\,|\,\filtration_t]$
are given by the relations
$$
(\widehat{\theta}_{1})_{{\rm BAYES}} =
        \frac{\Mprocess_t / \sigma + m / \gamma^2}{\wfunc(t) / \sigma^2 + 1/\gamma^2}
        \quad \mbox{è} \quad
{\bf E}\,\left[(\theta_1 - (\widehat{\theta}_{1})_{{\rm BAYES}})^2 \,|\,\filtration_t \right] =\frac{1}{\wfunc(t) / \sigma^2 + 1/\gamma^2}
$$
respectively, a result previously obrained by Norros~\cite{norros1999elementary} (for $\sigma = 1$).

\subsection{The case of a uniform prior distribution}\
Consider the problem of finding a Bayesian estimate 
for the value of the parameter~$\thetavector\in{\sR}^{n + 1}$
given that it is uniformly distributed on the~$(n + 1)$-dimensional cube
${\bf r} = \prod_{i=0}^n[a_i, b_i]$.

The density of the prior distribution of~$\thetavector$ is specified by
$$
p^{\thetavector}({\bf x})
    = \prod_{i=0}^{n} \frac{1}{b_i - a_i}\,{\bf 1}_{[a_i,b_i]}(x_i)
    =\frac{1}{|{\bf r}|} \,{\bf 1}_{{\bf r}}({\bf x}),
$$
where ${\bf 1}_{{\bf r}}({\bf x}) = \prod_{i=0}^{n} {\bf 1}_{[a_i,b_i]}(x_i)$,
$|{\bf r}|=\prod_{i=0}^{n} (b_i - a_i)^{-1}$. The corresponding posterior density
is given by the expression
$$
p^{\thetavector}({\bf x}\,|\,\filtration_t) =
\frac{1}{\Zprocess_t}\,{\bf 1}_{{\bf r}}({\bf x})
\exp\bigg\{
{\bf x}^{\top}\psiprocessn_t - \frac{1}{2}\,{\bf x}^{\top}\Rfunc(t){\bf x}\bigg\} = 
    \frac{1}{\Zprocess_t}\,{\bf 1}_{{\bf r}}({\bf x})
    \lprocess_t({\bf x}),
$$
where the process $\newprocess{\Zprocess}$, specified by the equality
$\Zprocess_t = \int_{{\bf r}} \lprocess_t({\bf x})\,d^n {\bf x}$,
plays the role of the normalization factor, and~$\lprocess_t({\bf x})$
is the likelihood process defined according to~\eqref{eq:fractional_likelihood_series}.

An analytic derivation of the normalization factor~$\Zprocess_t$,
the conditional mean~${\bf E}\,[\thetavector\,|\,\filtration_t]$
and the covariance matrix~${\rm cov}\,[\thetavector\,|\,\filtration_t]$
is a difficult problem for an arbitrary value of~$n$
(these quantities can be numerically computed using, e.\,g.,
the algorithm from~\cite{genz1992numerical}). We shall dwell upon
the derivation of the estimate for an important particular case
of a linear trend, where the observable process~$\xi$ satisfies
the following stochastic differential equation
\begin{equation}
d\xi_{t} = \theta_{1}\,dt + \sigma \,dB_{t}^{H},
\label{eq:uniform_bayes_n1_problem}
\end{equation}
where $\theta_1 \sim U(a, b)$. We present the result of our derivation
in the following theorem.

\begin{theorem}
\label{thm:bayesian_uniform}
Let~$\theta_1$ from~{\rm\eqref{eq:uniform_bayes_n1_problem}}
be a random variable uniformly distributed on $[a,b]$,
and independent of $B_{t}^{H}$. Then the optimal Bayesian
estimate of the value of the parameter~$\theta_1$ has the form
\begin{equation}
\label{eq:uniform_bayesian_estimate}
(\widehat{\theta}_{1})_{{\rm BAYES}} = \mprocess_t + [\Zprocess_t \wfunc(t)]^{-1}
[\lprocess_t(a) - \lprocess_t(b)],
\end{equation}
the conditional mean square estimation error is given by
\begin{eqnarray}
\label{eq:uniform_bayesian_estimate_var}
\gamma_t^H &=& {\bf E}\,\left((\theta_1 - (\widehat{\theta}_{1})_{{\rm BAYES}} )^2 \, |\, \filtration_t \right)\nonumber\\
&=& [\wfunc(t)]^{-1} + [Z_H (t)\wfunc(t)]^{-1}
[\lprocess_t(a)(a - \mprocess_t) - \lprocess_t(b)(b - \mprocess_t)]\nonumber\\
&&-\,[Z_H (t)\wfunc(t)]^{-2}
[ \lprocess_t(a) - \lprocess_t(b)]^2,
\end{eqnarray}
where
\begin{eqnarray}
\label{eq:uniform_bayesian_auxiliary_coefficients}
\Zprocess_t &=& \sqrt{\frac{2 \pi}{\wfunc(t)}}\,
\exp\bigg\{\frac{1}{2}\,(\mprocess_t\big)^2 \wfunc(t)\bigg\} \Cprocess_t,\\
\Cprocess_t &=& \Phi\left((b-\mprocess_t)\sqrt{\wfunc(t)}\,\right)
    - \Phi\left((a - \mprocess_t)\sqrt{\wfunc(t)}\,\right) \nonumber.
\end{eqnarray}
\end{theorem}

\textbf{Proof.} The conditional distribution~$p^{\theta_1}(x\,|\,\filtration_t)$
is easy to obtain using a direct computation, it is given by a formula
$$
p^{\theta_1}(x\,|\,\filtration_t)
 = \frac{1}{\Zprocess_t}\,{\bf 1}_{[a,b]}(x)
        \exp\bigg\{\wfunc(t) \Big(x \mprocess_t - \frac{x^2}{2}\Big)\bigg\},
$$
where the process~$\Zprocess$ is defined according
to~\eqref{eq:uniform_bayesian_auxiliary_coefficients}.
The conditional mean and variance are similarily obtained by
computing the corresponding integrals.

We further present several asymptotic properties of the obtained
Bayesian filter~\eqref{eq:uniform_bayesian_estimate}.

For $a \to -\infty, b \to +\infty$ (i.e.\ when $\theta_1$ is arbitrary)
the Bayesian estimate in~\eqref{eq:uniform_bayesian_estimate} coincides
with the maximum likelihood estimate. Indeed, as~$x \to \pm \infty$
we obtain~$\lprocess_t(x) \to 0$ so the second term
in~\eqref{eq:uniform_bayesian_estimate} vanishes
as $x \to \pm \infty$, meaning that~$(\widehat{\theta}_{1})_{{\rm BAYES}} \to m_t^H$.

As $t \to \infty$, the Bayesian estimate in~\eqref{eq:uniform_bayesian_estimate}
also coincides with the maximum likelihood estimate. Indeed,
for $t \to \infty$ we have~$\wfunc(t) \to \infty$, meaning that
the second term in~\eqref{eq:uniform_bayesian_estimate} vanishes
as~$t \to \infty$, and $(\widehat{\theta}_{1})_{{\rm BAYES}} \to m_t^H$.

Consider the problem of finding the optimal stopping time
in~\eqref{eq:optimality_criteria}. The cost function in this problem 
is given by
$$
{\bf E}\,\left[c\tau + {\bf E}\,\left((\theta_1
- (\widehat{\theta}_{1})_{{\rm BAYES}})^2 \, |\, \filtration_{\tau} \right) \right] = {\bf E\,}[c\tau + \gamma^H_{\tau}],
$$
where the random process $\newprocess{\gamma^H}$ is given by the
relation~\eqref{eq:uniform_bayesian_estimate_var}.
Note that as~$t \to \infty$, the following relation holds: $\gamma^H_t \to 0$.
To determine the optimal stopping time, it is neccessary to solve
\begin{equation}
\label{eq1}
\taubayes = \arg \inf_{\tau} {\bf E}\,[c\tau + \gamma^H_{\tau}].
\end{equation}

Since formulas~\eqref{eq:uniform_bayesian_estimate},
\eqref{eq:uniform_bayesian_estimate_var} and~\eqref{eq:uniform_bayesian_auxiliary_coefficients} 
are very complicated, an analytic solution for $\taubayes$ from~\eqref{eq1}
is infeasible, meaning that only a numerical estimation of the stopping
time is possible (see\ approaches in~\cite{shiryaev2006stop}).

\section{Acknowledgements}
The research was carried out at IITP RAS and supported by RSF grant No. 14-50-00150 only.

Artemov~A.\,V. is with Lomonosov Moscow State University,
Complex Systems Modelling Laboratory, Moscow, Russia;
Yandex Data Factory, Moscow, Russia (artemov@physics.msu.ru).

Burnaev~E.\,V. is with A.\,A.\,Kharkevich Institute of Information Transmission
Problems, RAS, Moscow, Russia (burnaev@iitp.ru).

\end{document}